\documentclass[12pt]{article}
\usepackage[fleqn]{amsmath}
\usepackage{amsfonts}

\newtheorem{theorem}{Theorem}[section]

\newtheorem{proposition}[theorem]{Proposition}

\newtheorem{remark}[theorem]{Remark}
\newtheorem{definition}[theorem]{Definition}

\DeclareMathOperator{\xGL}{GL}
\DeclareMathOperator{\xNS}{NS}
\DeclareMathOperator{\Tr}{Tr}

\DeclareMathOperator{\sgn}{sgn}

\begin{document}

\setlength{\pdfpagewidth}{8.5in}
\setlength{\pdfpageheight}{11in}

\title{Salem Numbers and Automorphisms of Abelian Surfaces}
\author{Paul Reschke\footnote{Research partially supported by NSF grants DMS0943832 and DMS1045119}}
\date{}

\maketitle

\begin{abstract}
We classify two-dimensional complex tori admitting automorphisms with positive entropy in terms of the entropies they exhibit. For each possible positive value of entropy, we describe the set of two-dimensional complex tori admitting automorphisms with that entropy.
\end{abstract}

\section{Overview}

Two-dimensional complex tori provide many basic examples of compact complex surfaces admitting (biholomorphic) automorphisms with positive topological entropy. Ghys and Verjovsky \cite{G_V} describe the circumstances under which a torus \(\mathbb{C}^2/\Lambda\) admits an infinite-order automorphism in terms of conditions on the lattice \(\Lambda\). Fujiki \cite{Fuj} describes the automorphism groups that can arise on two-dimensional complex tori, and indicates which types of tori allow which automorphism groups. Here, we characterize tori in terms of the entropies that they permit: for any fixed positive value that is the entropy of an automorphism on some two-dimensional complex torus, we describe the set of all such tori admitting automorphisms with the given entropy.

Suppose that \(f\) is an automorphism of a two-dimensional complex torus \(X=\mathbb{C}^2/\Lambda\), and let \(\gamma_1\) and \(\gamma_2\) be the eigenvalues of \(f^*\) on \(H^{1,0}(X) \cong \mathbb{C}^2\). So the eigenvalues of \(f^*\) on \(H^{0,1}(X)\) are \(\overline{\gamma_1}\) and \(\overline{\gamma_2}\). Since \(f^*\) is invertible on \(H^1(X,\mathbb{Z}) \cong \mathbb{Z}^4\), we must have \(|\gamma_1|^2 |\gamma_2|^2=1\). Moreover, since \(H^*(X,\mathbb{C})\) is generated by \(H^1(X,\mathbb{C})\) via the cup product, \(f^*\) has the following eigenvalues on \(H^2(X,\mathbb{C})\): \(\gamma_1 \gamma_2\) on \(H^{2,0}(X)\); \(\lambda = |\gamma_1|^2\), \(\gamma_1 \overline{\gamma_2}\), \(\overline{\gamma_1} \gamma_2\), and \(\lambda^{-1} = |\gamma_2|^2\) on \(H^{1,1}(X)\); and \(\overline{\gamma_1 \gamma_2}\) on \(H^{0,2}(X)\). Take \(\gamma_1\) and \(\gamma_2\) to be chosen so that \(\lambda \geq 1\). When \(f\) is algebraic (so \(f\) is the quotient of some \(F \in \xGL_2(\mathbb{C})\) with \(F(\Lambda)=\Lambda\)), it follows from a result of Sinai (\cite{Si1},\cite{Si2}) (applied to the four-dimensional real torus underlying \(X\)) that the topological entropy of \(f\) is \(\log(\lambda)\); in fact, the entropy is \(\log(\lambda)\) even if \(f\) is not algebraic. (See \S 2.2 below.) 

Since \(f^*\) has at most one eigenvalue with magnitude greater than one on \(H^2(X,\mathbb{Z})\), it follows from a result of Kronecker \cite{Kro} that the irreducible factors of the characteristic polynomial for \(f^*\) on \(H^2(X,\mathbb{Z})\) can only be cyclotomic polynomials and at most one Salem polynomial (an irreducible polynomial that is monic and reciprocal and has exactly two roots with magnitude not equal to one); thus \(f\) has positive entropy if and only if \(\lambda\) is a Salem number (the largest real root of a Salem polynomial). In \cite{Re1}, we characterized all of the Salem polynomials that give entropies of two-dimensional complex torus automorphisms. This paper completes that effort by describing the tori on which these Salem polynomials arise. Note that the degree of any such polynomial must be two, four, or six.

\begin{theorem}
Suppose \(\lambda\) is a Salem number such that \(\log(\lambda)\) is the entropy of some two-dimensional complex torus automorphism, and let \(S(t)\) be the minimal polynomial for \(\lambda\).
\begin{itemize}
\item[1)] If \(\deg(S(t))=6\), then any two-dimensional complex torus automorphism on which \(\log(\lambda)\) is the entropy of an automorphism must be non-projective.
\item[2)] If \(\deg(S(t))=4\), then \(\log(\lambda)\) is the entropy of an automorphism of an abelian surface and the entropy of an automorphism of a non-projective two-dimensional complex torus.
\item[3)] If \(\deg(S(t))=2\), then any two-dimensional complex torus automorphism on which \(\log(\lambda)\) is the entropy of an automorphism must be projective.
\begin{itemize}
\item[a)] If either \(\lambda + \lambda^{-1} + 2\) or \(\lambda + \lambda^{-1} - 2\) is the square of an integer, then \(\log(\lambda)\) is the entropy of an automorphism of a simple abelian surface and the entropy of an automorphism of a product of isogenous elliptic curves.
\item[b)] Otherwise, any abelian surface on which \(\log(\lambda)\) is the entropy of an automorphism must be a product of isogenous elliptic curves with complex multiplication.
\end{itemize}
\end{itemize}
\end{theorem}
\noindent (See \S 3.2 below for the proof.)

Note that (1) actually follows from a straightforward observation: for an automorphism \(f\) of a two-dimensional complex torus \(X\) with \(\xNS(X) \neq \{0\}\), the characteristic polynomial for \(f^*\) on \(\xNS(X) \subseteq H^{1,1}(X)\) must be a factor (of degree at most four) in the characteristic polynomial for \(f^*\) on \(H^2(X,\mathbb{Z})\). Indeed, every torus in (1) must have Picard rank zero. On the other hand, every non-projective torus in (2) has Picard rank two. These tori turn out to be generalizations of examples developed by Zucker \cite{Zuc} to show that a non-projective compact K\"ahler manifold can have a non-trivial N\'eron-Severi group without having any divisors. (See \S 4.3 below.) Every abelian surface in (2) has Picard rank four; by work of Shioda and Mitani \cite{S_M}, it follows that each such torus is a product of isogenous elliptic curves with complex multiplication. Every abelian surface in (3) is either simple or isogenous to a product of isogenous elliptic curves (possibly without complex multiplication). (See \S 4.1 below.)

The set of entropies arising on a torus of the form \(E \times E\), where \(E\) is an elliptic curve, is complicated--especially when \(E\) has complex multiplication; on the other hand, the set of entropies exhibited by the automorphism group of a non-projective two-dimensional complex torus is either trivial or equal to
\[\{k\log(\lambda)|k \in \mathbb{N}_0\}\]
for some Salem number \(\lambda\). (See \S 5.1 below.) In the opposite direction, we show that the entropy of a two-dimensional complex torus automorphism will typically only occur on finitely many two-dimensional complex tori.

\begin{theorem}
Let \(\lambda\) be a Salem number such that \(\log(\lambda)\) is the entropy of some two-dimensional complex torus automorphism. Then either:
\begin{itemize}
\item[1)] One of \(\lambda+\lambda^{-1}+2\) or \(\lambda+\lambda^{-1}-2\) is the square of an integer; or
\item[2)] The set of two-dimensional complex tori that admit automorphisms whose entropies are \(\log(\lambda)\) is finite.
\end{itemize}
\end{theorem}
\noindent (See \S 5.3 below for the proof.)

For each \(\lambda\) in (1), there is a positive-dimensional set of parameters defining abelian surfaces--including both simple and non-simple abelian surfaces--that admit automorphisms whose entropies are \(\log(\lambda)\); indeed, there is a Hilbert modular surface which parametrizes a subset of all such abelian surfaces. (See \S 4.1 below.) The finite number of isomorphism classes of tori in (2) can be arbitrarily large as \(\lambda\) varies; however, the number of isogeny classes of tori in (2) is uniformly bounded for all \(\lambda\). (See \S 5.3 below.)
\\
\\
\noindent \textbf{Acknowledgements.} Curt McMullen provided some of the initial motivation for this paper (by asking questions about a previous result of the author in \cite{Re1}) and subsequently made several helpful suggestions and comments. Christian Schnell provided useful details regarding Zucker's examples of non-projective complex tori with non-trivial N\'eron-Severi groups and no divisors. Igor Dolgachev pointed out the relevance of Hilbert modular surfaces to Theorem 1.2. Holly Krieger performed calculations for a preliminary result which was eventually superseded by Theorem 1.1.

\section{Topological Entropy on Compact K\"ahler Surfaces}

The following theorem is a special case of a compilation of results by Gromov \cite{Gro}, Yomdin \cite{Yom}, and Friedland \cite{Fri}. (See also \cite{Gue}, \S 2.)

\begin{theorem}
Let \(f\) be an automorphism of a compact K\"ahler surface \(X\). Then the topological entropy of \(f\) is the logarithm spectral radius of \(f^*\) on \(H^{1,1}(X)\).
\end{theorem}

The conclusion that the entropy of \(f\) is either zero or the logarithm of a Salem number is not special to the case where \(X\) is a torus; indeed, the irreducible factors of the characteristic polynomial for \(f^*\) on \(H^2(X,\mathbb{Z})\) can only be cyclic polynomials and at most one Salem polynomial for any \(f\) and \(X\) in Theorem 2.1. (See, e.g., \cite{Mc1}, \S 3.)

\subsection{Positive Entropy on Compact K\"ahler Surfaces}

Given an automorphism \(f\) of a compact K\"ahler surface \(X\), we can use the following theorem to determine whether \(X\) is projective or not based solely on the cohomological behavior of \(f\); the statement that (2) implies (1) is new, while the statement that (1) implies (2) follows immediately from previously known results.

\begin{theorem}
Let \(f\) be an automorphism of a compact K\"ahler surface \(X\). Suppose that the entropy of \(f\) is \(\log(\lambda)>1\), and let \(S(t)\) be the minimal polynomial for \(\lambda\). Then the following two statements are equivalent:
\begin{itemize}
\item[1)] \(H^{2,0}(X) \neq \{0\}\) and some root of \(S(t)\) is an eigenvalue for \(f^*\) on \(H^{2,0}(X)\); and
\item[2)] \(X\) is non-projective.
\end{itemize}
\end{theorem}

\emph{Proof:} Each root of \(S(t)\) is a simple eigenvalue for \(f^*\) on \(H^2(X,\mathbb{Z})\). If \(X\) is projective, then \(\lambda\) is an eigenvalue for \(f^*\) on \(\xNS(X)\); but then every root of \(S(t)\) is an eigenvalue for \(f^*\) on \(\xNS(X)\) (and hence not an eigenvalue for \(f^*\) on \(H^{2,0}(X)\)). (See \cite{Re2}, \S 3.2, and \cite{Mc1}, \S 3.) So (1) implies (2).

Since \(S(t)\) is the only non-cyclotomic factor for \(f^*\) on \(H^2(X,\mathbb{Z})\), there is a subspace \(W \subseteq H^2(X,\mathbb{Q})\) such that \(f^*(W)=W\) and the characteristic polynomial for \(f^*\) on \(W\) is \(S(t)\). (See \cite{Re2}, \S 3.2.) So \(W \otimes \mathbb{C}\) is the span of the set of all eigenvectors for \(f^*\) on \(H^2(X,\mathbb{C})\) corresponding to eigenvalues which are roots of \(S(t)\). If every root of \(S(t)\) is an eigenvalue for \(f^*\) on \(H^{1,1}(X)\), then \(W \otimes \mathbb{C}\) is a subspace of \(H^{1,1}(X)\); it follows from the Lefschetz theorem on \((1,1)\)-classes that \(W\) is a subspace of \(\xNS(X) \otimes \mathbb{Q}\). Thus, since \(f^*\) preserves the intersection pairing on \(\xNS(X)\), \(\xNS(X) \otimes \mathbb{R}\) must contain two distinct totally isotropic subspaces of dimension one (consisting of the eigenvectors corresponding to the eigenvalues \(\lambda\) and \(\lambda^{-1}\)); since the signature of the real part of \(H^{1,1}(X)\) is \((1,h^{1,1}(X)-1)\), it follows that \(\xNS(X)\) must contain some element with positive self-intersection and hence that \(X\) must be projective. (See \cite{Re2}, \S 3.1.) So (2) implies (1). \(\Box\)

\begin{remark}
Cantat \cite{Ca1} showed that any compact K\"ahler surface admitting an automorphism with positive entropy must be birational to a torus, a K3 surface, an Enriques surface, or \(\mathbb{P}^2\). This result limits the scope of Theorem 2.2 to (blow-ups of) tori and K3 surfaces. However, the proof of Theorem 2.2 does not rely a priori on any constraint on \(X\).
\end{remark}

\subsection{Positive Entropy on Two-Dimensional Complex Tori}

The following proposition shows that we may restrict our attention to algebraic automorphisms when considering two-dimensional complex torus automorphisms with positive entropy. An automorphism of a torus \(X=\mathbb{C}^2/\Lambda\) is algebraic if it is the quotient of some \(F \in \xGL_2(\mathbb{C})\) satisfying \(F(\Lambda)=\Lambda\)--or, equivalently, if it respects the group structure on \(X\) inherited from \(\mathbb{C}^2\).

\begin{proposition}
Let \(f\) be an automorphism of a two-dimensional complex torus \(X=\mathbb{C}^2/\Lambda\), let \(\gamma_1\) and \(\gamma_2\) be the eigenvalues of \(f^*\) on \(H^{1,0}(X)\), and suppose that \(|\gamma_1| > 1\). Then \(f\) is conjugate by a translation on \(X\) to an algebraic automorphism of \(X\).
\end{proposition}

\emph{Proof:} There is an algebraic automorphism \(\phi_f:X \rightarrow X\) and an element \(x_f \in X\) such that
\[f(x)=\phi_f(x)+x_f\]
for any \(x \in X\). (See, e.g., \cite{B_L}, \S 1.1.)  The Lefschetz number for \(f\) is
\[\sum (-1)^j \Tr(f^*:H^j(X,\mathbb{Z}) \rightarrow H^j(X,\mathbb{Z})) = (1-\gamma_1)(1-\gamma_2)(1-\overline{\gamma_1})(1-\overline{\gamma_2}) \neq 0.\]
Thus the Lefschetz fixed-point theorem guarantees that \(f(x_0)=x_0\) for some \(x_0 \in X\). (See also \cite{Zha}, \S 2.1.) So
\[f(x+x_0)-x_0 = \phi_f(x)+\phi_f(x_0)+x_f-x_0 = \phi_f(x)\]
for any \(x \in X\); that is, \(f\) is conjugate (by translation by \(x_0\)) to \(\phi_f\). \(\Box\)

\begin{remark}
Berg \cite{Ber} showed that Haar measure is always a measure of maximal entropy for an algebraic automorphism of a real torus; thus the work of Sinai (\cite{Si1},\cite{Si2}) on entropies of algebraic automorphisms with respect to Haar measure gives the conclusion of Theorem 2.1 for algebraic automorphisms of two-dimensional complex tori. Proposition 2.4 then gives this conclusion for any two-dimensional complex torus automorphism with cohomological eigenvalues which are not roots of unity. So, for the proofs of Theorems 1.1 and 1.2, the general statement of Theorem 2.1 is necessary only for the conclusion that any two-dimensional complex torus automorphism (algebraic or not) whose cohomological eigenvalues are all roots of unity has entropy zero.
\end{remark}

\section{Types of Tori Exhibiting Positive Entropies}

Suppose that \(f\) is an algebraic automorphism with positive entropy of a torus \(X = \mathbb{C}^2/\Lambda\); so \(f\) is the quotient of some \(F \in \xGL_2(\mathbb{C})\) satisfying \(F(\Lambda)=\Lambda\). Since \(H^{1,0}(X)\) is spanned by \(dz_1\) and \(dz_2\) (for any choice of coordinate system \(\{z_1,z_2\}\) on \(\mathbb{C}^2\)), it follows that \(f^* = F^T\) on \(H^{1,0}(X)\). 

\subsection{Reorientations of Tori by Automorphisms}

\begin{definition}
Suppose that \(F \in \xGL_2(\mathbb{C})\) satisfies \(F(\Lambda)=\Lambda\) for some lattice \(\Lambda \subseteq \mathbb{C}\) and has eigenvalues with magnitude different from one, and let \(f\) be the automorphism (with positive entropy) of \(X = \mathbb{C}^2/\Lambda\) which is the quotient of \(F\). Choose a basis for \(\mathbb{C}^2\) with respect to which
\[F= \left(
\begin{array}{cc}
\gamma_1 & 0 \\
0 & \gamma_2
\end{array} \right)\]
(where \(\gamma_1\) and \(\gamma_2\) are the eigenvalues of \(F\)), and set
\[\Lambda' = \{(z_1,z_2)|(z_1,\overline{z_2}) \in \Lambda\}\]
in this basis. Then
\[F'= \left(
\begin{array}{cc}
\gamma_1 & 0 \\
0 & \overline{\gamma_2}
\end{array} \right)\]
satisfies \(F'(\Lambda')=\Lambda'\), and hence induces an automorphism \(f'\) of \(X' = \mathbb{C}^2/\Lambda'\). We say that \(f'\) is the \emph{reorientation} of \(f\), and that \(X'\) is the \emph{reorientation} of \(X\) by \(f\).
\end{definition}

In Definition 3.1, the reorientations \(f'\) and \(X'\) are independent of the choice of basis diagonalizing \(F\). Note that \(f\) and \(f'\) have the same entropy.

\begin{proposition}
Suppose that \(f\) is an algebraic automorphism of a two-dimensional complex torus \(X\) whose entropy is the logarithm of a degree-four Salem number \(\lambda\), and let \(X'\) be the reorientation of \(X\) by \(f\). Then exactly one of \(X\) or \(X'\) is projective.
\end{proposition}

\emph{Proof:} Let \(S(t)\) be the minimal polynomial for \(\lambda\); so the characteristic polynomial for \(f^*\) on \(H^2(X,\mathbb{Z})\) has the form
\[S(t)(t^2+at+1)\]
for some \(a \in \{-2,-1,0,1,2\}\). Let \(\gamma_1\) and \(\gamma_1\) be the eigenvalues for \(f^*\) on \(H^{1,0}(X)\); then exactly one of \(\gamma_1 \gamma_2\) or \(\gamma_1 \overline{\gamma_2}\) is a root of unity. Since \(\gamma_1 \gamma_2\) is the eigenvalue for \(f^*\) on \(H^{2,0}(X)\) and \(\gamma_1 \overline{\gamma_2}\) is the eigenvalue for \((f')^*\) on \(H^{2,0}(X')\) (where \(f'\) is the reorientation of \(f\)), it follows from Theorem 2.2 that exactly one of \(X\) of \(X'\) is projective. \(\Box\)

\subsection{Proof of Theorem 1.1}

\noindent \emph{Proof of Theorem 1.1:} If a two-dimensional complex torus admits an automorphism whose entropy is \(\log(\lambda)\), then it follows from Proposition 2.4 that the torus also admits an algebraic automorphism whose entropy is \(\log(\lambda)\). Suppose that \(X\) is one such torus, and that \(f\) is one such algebraic automorphism. Since the real part of
\[H^{2,0}(X) \oplus H^{0,2}(X)\]
has signature \((2,0)\), the eigenvalue for \(f^*\) on \(H^{2,0}(X)\) must have magnitude one. (See, e.g., \cite{Mc1}, \S 3.) Let \(Q(t)\) and \(P(t)\) be the characteristic polynomials for \(f^*\) on \(H^2(X,\mathbb{Z})\) and \(H^1(X,\mathbb{Z})\), respectively. Since the roots of \(Q(t)\) are precisely the products of distinct pairs of roots of \(P(t)\), we can compute the coefficients of \(Q(t)\) in terms of the coefficients of \(P(t)\) and observe that there are integers \(m\) and \(n\) such that \(Q(1)=-m^2\) and \(Q(-1)=n^2\); we have also that \(Q(t)\) is monic and reciprocal. (See \cite{Re1}, \S 3.)

If \(\deg(S(t))=6\), then every eigenvalue for \(f^*\) on \(H^2(X,\mathbb{Z})\) is a root of \(S(t)\); so the conclusion follows from Theorem 2.2. If \(\deg(S(t))=4\), then the conclusion follows from Proposition 3.2.

If \(\deg(S(t))=2\), then the eigenvalues for \(f^*\) on \(H^2(X,\mathbb{Z})\) are \(\lambda\), \(\lambda^{-1}\), and four roots of unity; so it follows from Theorem 2.2 that \(X\) is an abelian surface. In this case, \(Q(t)=S(t)C(t)\) for some monic and reciprocal degree-four polynomial \(C(t)\) whose roots are all roots of unity; specifically, \(C(t)\) must either factor as
\[(t^2+jt+1)(t^2+kt+1)\]
for some \(j\) and \(k\) in \(\{-2,-1,0,1,2\}\) or be in
\[\{t^4+t^3+t^2+t+1,t^4+1,t^4-t^3+t^2-t+1,t^4-t^2+1\}.\]
Suppose that neither
\[-S(1) = \lambda + \lambda^{-1}-2\]
nor
\[S(-1) = \lambda + \lambda^{-1}+2\]
is the square of an integer. Then the conditions \(Q(1)=-m^2\) and \(Q(-1)=n^2\) imply that \(C(t)\) must be reducible with \(j \neq k\); moreover, since the eigenvalues for \(f^*\) on \(H^{2,0}(X)\) and \(H^{0,2}(X)\) are complex conjugates of one another, \(j\) and \(k\) can be ordered so that the characteristic polynomial for \(f^*\) on \(H^{1,1}(X)\) is
\[Q_0(t) = S(t)(t^2+kt+1)\]
(and the characteristic polynomial for \(f^*\) on \(H^{2,0}(X) \oplus H^{0,2}(X)\) is \(t^2+jt+1\)).
As in the proof of Theorem 2.2, there is a subspace \(W \subseteq H^2(X,\mathbb{Q})\) such that \(f^*(W)=W\) and the characteristic polynomial for \(f^*\) on \(W\) is \(Q_0(t)\); it follows from the Lefschetz theorem on (1,1)-classes that \(X\) must have Picard rank four (since \(W \otimes \mathbb{C} = H^{1,1}(X)\)). So \(X\) must be a product of isogenous elliptic curves with complex multiplication if the hypothesis of (a) fails. (See \cite{S_M}, \S 4.) Example 4.1 below demonstrates the conclusion of (a) if its hypothesis holds. \(\Box\)

\begin{remark}
Since reorientation does not change entropy, Theorem 1.1 shows that the conclusion of Proposition 3.2 is special to the case of degree-four Salem numbers: reorientation by an automorphism whose entropy is the logarithm of a degree-two Salem number must interchange two abelian surfaces, while reorientation by an automorphism whose entropy is the logarithm of a degree-six Salem number must interchange two non-projective tori.
\end{remark}

\section{Examples of Tori Exhibiting Positive Entropies}

\subsection{Case 3a in Theorem 1.1}

Let \(\lambda\) be a degree-two Salem number; so
\[q = \lambda+\lambda^{-1}\]
is an integer greater than \(2\), and the minimal polynomial for \(\lambda\) is
\[S(t)=t^2-qt+1.\]
Suppose that \(q+2=r^2\) (resp., \(q-2=r^2\)) for some integer \(r\), and let \(A =(a_{ij}) \in \xGL_2(\mathbb{Z})\) be a matrix with determinant \(1\) (resp., \(-1\)) and trace \(r\); so the eigenvalues of \(A\) are \(\sgn(r)\sqrt{\lambda}\) and \(\sgn(r)\sqrt{\lambda}^{-1}\) (resp., \(\sgn(r)\sqrt{\lambda}\) and \(-\sgn(r)\sqrt{\lambda}^{-1}\)). Then any two-dimensional complex torus of the form \(E \times E\), where \(E\) is an elliptic curve, admits an automorphism \(\sigma\) given by
\[\sigma(e_1,e_2) = (a_{11}e_1+a_{12}e_2,a_{21}e_1+a_{22}e_2)\]
whose entropy is \(\log(\lambda)\). (See also \cite{Mc1}, \S 4, and \cite{Kaw}, \S 3.)

By the Poincar\'e reducibility theorem, every abelian surface is either simple or isogenous to \(E_1 \times E_2\) for some elliptic curves \(E_1\) and \(E_2\). If \(E_1\) and \(E_2\) are two non-isogenous elliptic curves, then any algebraic endomorphism of \(E_1 \times E_2\) is given by
\[(e_1,e_2)=(f_1(e_1),f_2(e_2))\]
for some algebraic endomorphisms \(f_1\) and \(f_2\) of, respectively, \(E_1\) and \(E_2\); if \(A\) is an abelian surface isogenous to \(E_1 \times E_2\), it follows that any algebraic automorphism of \(A\) must leave invariant the images of both \(E_1 \times \{0\}\) and \(\{0\} \times E_2\)--and hence must have eigenvalues on \(H^{1,0}(A)\) which are roots of unity (since every algebraic automorphism of an elliptic curve has finite order). (See, e.g, \cite{Mum}, \S IV.19, and \cite{Mir}, \S III.1.) So every abelian surface that admits an automorphism with entropy \(\log(\lambda)\) is either simple or isogenous to some \(E \times E\) (with \(E\) an elliptic curve).

Let \(\gamma_1\) and \(\gamma_2\) be the eigenvalues of \(A\) (so \(r=\gamma_1+\gamma_2\) and \(\det(A)=\gamma_1 \gamma_2\)), and suppose that \(\Lambda\) is a lattice in \(\mathbb{C}^2\) with a basis of the form
\[\{(1,1),(\gamma_1,\gamma_2),(z_1,z_2),(\gamma_1 z_1,\gamma_2 z_2)\}\]
(for some complex numbers \(z_1\) and \(z_2\)). Since \(\gamma_1\) and \(\gamma_2\) are both roots of \(t^2-rt+\det(A)\),
\[\left( \begin{array}{cc}
\gamma_1 & 0 \\
0 & \gamma_2
\end{array} \right)\]
\noindent
restricts to
\[\left( \begin{array}{cccc}
0 & -\det(A) & 0 & 0 \\
1 & r & 0 & 0 \\
0 & 0 & 0 & -\det(A) \\
0 & 0 & 1 & r
\end{array} \right)\]
\noindent
on \(\Lambda\); so \(\mathbb{C}^2/\Lambda\) admits an automorphism whose eigenvalues on \(H^{1,0}(\mathbb{C}^2/\Lambda)\) are \(\gamma_1\) and \(\gamma_2\) (and whose entropy is therefore \(\log(\lambda)\)). Varying \(z_1\) and \(z_2\) yields a positive-dimensional set of parameters describing abelian surfaces that admit automorphisms whose entropries are \(\log(\lambda)\).

Suppose further that \(z_1=\imath\) and \(z_2=\delta \imath\) for some non-zero \(\delta \in \mathbb{R}\). The abelian surface \(\mathbb{C}^2/\Lambda\) contains an elliptic curve if and only if it contains two distinct isogenous elliptic curves, in which case there are elements
\[(\zeta_1,\zeta_2) = (k_1,k_1)+(k_2 \gamma_1,k_2 \gamma_2)+(k_3 \imath,k_3 \delta \imath) \in \Lambda\]
(so each \(k_j\) is an integer) and \(c+d\imath \in \mathbb{C}\) with \(d \neq 0\) such that \((c+d\imath)(\zeta_1,\zeta_2) \in \Lambda\). If the equations
\[ck_1+ck_2\gamma_1-dk_3+dk_1\imath+dk_2\gamma_1\imath+ck_3\imath = l_1+l_2\gamma_1+l_3\imath+l_4\gamma_1\imath\]
and
\[ck_1+ck_2\gamma_2-dk_3\delta +dk_1\imath+dk_2\gamma_2\imath+ck_3\delta\imath = l_1+l_2\gamma_2+l_3\delta\imath+l_4\gamma_2\delta\imath\]
have a simultaneous non-trivial solution with \(k_1,\dots,k_3,l_1,\dots,l_4 \in \mathbb{Z}\) and \(c,d \in \mathbb{R}\) (with \(d \neq 0\)), then
\[dk_3(\delta k_1+\delta k_2\gamma_1-k_1-k_2\gamma_2) = (l_2k_1-l_1k_2)(\gamma_1-\gamma_2) = l_4k_3\delta(\gamma_1-\gamma_2)\]
implies \(\delta \in \mathbb{Q}\) or
\[\delta(k_1+k_2\gamma_1)=k_1+k_2\gamma_2\]
or
\[\delta(l_3+l_4\gamma_2)(k_1+k_2\gamma_1)=(l_3+l_4\gamma_1)(k_1+k_2\gamma_2).\]
So \(\mathbb{C}^2/\Lambda\) is simple for a generic choice of \(\delta\) (including, for example, any \(\delta\) not contained in \(\mathbb{Q}(\gamma_1)\)).

Since
\[\lambda = (q+\sqrt{q^2-4})/2,\]
the condition \(q \pm 2 = r^2\) is equivalent to the condition that \(\sqrt{\lambda}\) is again a quadratic integer; the minimal polynomial for \(\sqrt{\lambda}\) is
\[t^2-rt \pm 1,\]
so that \(\lambda\) is Galois conjugate to \(\pm \sqrt{\lambda}^{-1}\). Thus
\[\mathbb{Q}(\sqrt{\lambda}) = \mathbb{Q}(\lambda) = \mathbb{Q}(\sqrt{q^2-4}) = \mathbb{Q}(\sqrt{D})\]
for some square-free \(D \in \mathbb{N}\), and \(\sqrt{\lambda}\) is a unit in the ring of integers for \(\mathbb{Q}(\sqrt{D})\). Suppose that \(X = \mathbb{C}^2/\Lambda\) is an abelian surface with multiplication by \(\sqrt{D}\); so there is a basis for \(\mathbb{C}^2\) such that
\[F_{a,b} = \left( \begin{array}{cc}
a+b\sqrt{D} & 0 \\
0 & a-b\sqrt{D}
\end{array} \right)\]
satisfies \(F_{a,b}(\Lambda) \subseteq \Lambda\) for any integer \(a+b\sqrt{D} \in \mathbb{Q}(\sqrt{D})\). (See \cite{V_U}, \S 1.) Taking \(\sqrt{\lambda} = a+b\sqrt{D}\) then gives an automorphism of \(X\) whose entropy is \(\log(\lambda)\). So the set of all abelian surfaces with multiplication by \(\sqrt{D}\)--which is parametrized by a Hilbert modular surface--is a subset of the set of all abelian surfaces admitting automorphisms with entropy \(\log(\lambda)\). On the other hand, if
\[F = \left( \begin{array}{cc}
\sqrt{\lambda} & 0 \\
0 & \pm \sqrt{\lambda}^{-1}
\end{array} \right) \in \xGL_2(\mathbb{C})\]
satisfies \(F(\Lambda)=\Lambda\) for some lattice \(\Lambda \subseteq \mathbb{C}^2\), then there is some integer \(b\) such that
\[F' = \left( \begin{array}{cc}
b\sqrt{D} & 0 \\
0 & -b\sqrt{D}
\end{array} \right)\]
satisfies \(F'(\Lambda) \subseteq \Lambda\). So any abelian surface admitting an automorphism with entropy \(\log(\lambda)\) must have real multiplication by \(b\sqrt{D}\).

\subsection{Case 3b in Theorem 1.1}

Let \(\lambda\) be a degree-two Salem number that does not satisfy the hypothesis of Example 4.1, and let \(S(t)\) be the minimal polynomial for \(\lambda\). Then the proof of Theorem 1.1 shows that the characteristic polynomial for \(f^*\) on \(H^2(X,\mathbb{Z})\) is
\[S(t)(t^2+jt+1)(t^2+kt+1)\]
for some distinct \(j\) and \(k\) in \(\{-2,-1,0,1,2\}\). In fact, taking \(j=-2\) and \(k=2\) will always give a polynomial that occurs as the characteristic polynomial on the second cohomology group for some abelian surface automorphism with positive entropy. (See \cite{Re1}, \S 3.)

Let \(E=\mathbb{C}/\mathbb{Z}[\sqrt{2-q}]\) (where \(q=\lambda+\lambda^{-1}\), as in Example 4.1). Then \(E\) has complex multiplication (by \(\sqrt{2-q}\)) and
\[\left( \begin{array}{cc}
0 & -1 \\
1 & (\sqrt{q-2})\imath
\end{array} \right)\]
\noindent
gives an automorphism of \(E \times E\) (as in Example 4.1) whose entropy is \(\log(\lambda)\).

\subsection{Case 2 in Theorem 1.1}

If \(X\) is an abelian surface that admits an automorphism whose entropy is the logarithm of a degree-four Salem number, then the proof of Theorem 2.2 shows that the Picard rank of \(X\) is four. If \(X\) is a non-projective two-dimensional complex torus that admits an automorphism \(\sigma\) whose entropy is the logarithm of a degree-four Salem number \(\lambda\), so that the characteristic polynomial for \(\sigma^*\) on \(H^2(X,\mathbb{Z})\) has the form
\[S(t)(t^2 + a t + 1)\]
with \(S(\lambda)=0\) and \(a \in \{-2,-1,0,1,2\}\), then Theorem 2.2 shows that \(t^2 + a t + 1\) is a factor in the characteristic polynomial for \(\sigma^*\) on \(H^{1,1}(X)\); it follows that the Picard rank of \(X\) is two.

Let \(E=\mathbb{C}/\mathbb{Z}[\sqrt{-D}]\) with \(D \in \mathbb{N}\). Then any matrix of the form
\[\left( \begin{array}{cc}
0 & -1 \\
1 & b_1 + (b_2 \sqrt{D})\imath
\end{array} \right),\]
\noindent
where \(b_1\) and \(b_2\) are integers, gives an automorphism of \(E \times E\) (as in Examples 4.1 and 4.2); the characteristic polynomial for the action of the automorphism on \(\xNS(E \times E)\) is
\[t^4 - (b_1^2 + b_2^2 D)t^3 +(2b_1^2 - 2b_2^2 D - 2)t^2 - (b_1^2 + b_2^2 D)t + 1,\]
which is a degree-four Salem polynomial whenever it does not have the form
\[t^4-pt^3-(2 \pm 2p)t^2-pt+1 \text{, } t^4-pt^3+(1 \pm p)t^2-pt+1 \text{, or } t^4-pt^3+2t^2-pt+1\]
(for \(p \in \mathbb{N}_0\)). (See \cite{Bea}, \S 5.2.) Let \(A_{b_1,b_2}\) be such a matrix, and let \(\sigma_{b_1,b_2}\) be the corresponding automorphism of \(E \times E\). Then the eigenvectors of \(A_{b_1,b_2}\) are
\[\left( 1,\frac{-b_1-b_2\imath \pm \sqrt{b_1^2-b_2^2+2 b_1 b_2 \imath -4}}{2} \right),\]
and the reorientation of \(E \times E\) by \(\sigma_{b_1,b_2}\) can be given concretely via an explicit change of basis for \(A_{b_1,b_2}\); since any lattice giving \(E \times E\) as a quotient of \(\mathbb{C}^2\) is invariant under the map that sends \((z_1,z_2)\) to \((\sqrt{D}\imath z_1,\sqrt{D}\imath z_2)\), the lattice giving the reorientation must be invariant under the map that sends \((z_1,z_2)\) to \((\sqrt{D}\imath z_1,-\sqrt{D}\imath z_2)\).

If \(\Lambda \subseteq \mathbb{C}^2\) is a lattice that is invariant under the map that sends \((z_1,z_2)\) to \((\sqrt{D}\imath z_1,\linebreak -\sqrt{D}\imath z_2)\) (with \(D \in \mathbb{N}\)), so that \(\Lambda\) has a basis of the form
\[\{(u_1,u_2),(v_1,v_2),(\sqrt{D}\imath u_1,-\sqrt{D}\imath u_2),(\sqrt{D}\imath v_1,-\sqrt{D}\imath v_2)\},\]
then there is a form \(\omega\) on \(\mathbb{C}^2/\Lambda\) such that \([\Re \omega]\) and \([\sqrt{D}\Im \omega]\) span a two-dimensional lattice in\ \(\xNS(X)\); in terms of the chosen basis for \(\mathbb{C}^2\),
\[\omega = (u_1 \overline{v_2} - v_1 \overline{u_2})^{-1} dz_1 \wedge d\overline{z_2}\]
has this property. A torus that can be expressed as the quotient of \(\mathbb{C}^2\) by a lattice that is invariant under the map that sends \((z_1,z_2)\) to \((\sqrt{D}\imath z_1,-\sqrt{D}\imath z_2)\) is called a \(J_D\)-torus; the fact that (for any \(D \in \mathbb{N}\)) certain \(J_D\)-tori admit automorphisms whose entropies are logarithms of degree-four Salem numbers shows that a generic \(J_D\) torus has Picard rank two. The intersection form is negative definite on \(\xNS(X)\) for any \(J_D\)-torus \(X\) with Picard rank two; so, since a two-dimensional complex torus cannot contain a curve with negative self-intersection (because the adjunction formula implies that such a curve would necessarily be an embedding of \(\mathbb{P}^1\) into the torus), a generic \(J_D\)-torus (including any \(J_D\)-torus that admits an automorphism with positive entropy) has no divisors.

\begin{remark}
Zucker (\cite{Zuc}, Appendix B) showed that a generic \(J_1\)-torus has Picard rank two and no divisors. Thus one application of the idea of reorientation by an automorphism with positive entropy is to give an alternate proof of Zucker's result.
\end{remark}

\subsection{Case 1 in Theorem 1.1}

For any \(a \in \mathbb{N}_0\),
\[S(t) = t^6-at^5-t^4+(2a-1)t^3-t^2-at+1\]
is a Salem polynomial whose roots are precisely the products of distinct pairs of roots of
\[P(t) = t^4+at^2+t+1;\]
it follows that there is a two-dimensional complex torus with an automorphism whose entropy is the logarithm of the Salem root of \(S(t)\). (See \cite{Mc1}, \S 4.) Indeed, any matrix of the form
\[A_b = \left( \begin{array}{cccc}
0 & 0 & -1 & 0 \\
1 & 0 & 0 & 0 \\
0 & 1 & 1 & b\\
0 & 0 & -(1+a)/b & -1
\end{array} \right),\]
where \(b\) is an integer dividing \(1+a\), is an invertible transformation (among others) of \(\mathbb{Z}^4\) with characteristic polynomial \(P(t)\); since the roots of \(P(t)\) are complex (and hence occur in complex conjugate pairs), \(\mathbb{Z}^4 \otimes \mathbb{R}\) can be given a complex structure (isomorphic to \(\mathbb{C}^2\)) on which \(A_b\) acts as an element of \(\xGL_2(\mathbb{C})\) preserving the lattice which is the image of \(\mathbb{Z}^4\). A similar construction applies for any degree-six Salem number whose logarithm is the entropy of some two-dimensional complex torus automorphism. (See \cite{Re2}, \S 3.)

\section{Limitations on Occurences of Positive Entropies}

\subsection{The Set of Entropies on a Fixed Torus}

For a non-projective two-dimensional complex torus \(X\) with an infinite-order automorphism, Fujiki (\cite{Fuj}, \S 5) asserts (without proof) that the algebraic automorphism group of \(X\) is isomorphic to
\[\mathbb{Z} \times (\mathbb{Z}/m\mathbb{Z})\]
for some \(m \in \{2,4,6\}\)--so that, in particular, there is some (infinite-order) automorphism \(\sigma_0\) of \(X\) such that the algebraic automorphism group of \(X\) is generated by \(\sigma_0\) and finitely many finite-order automorphisms; since any two algebraic automorphisms of \(X\) commute with one another in this case, it follows that the set of entropies exhibited by the automorphism group of \(X\) is
\[\{k\log(\lambda_0)|k \in \mathbb{N}_0\},\]
where \(\log(\lambda_0)\) is the entropy of \(\sigma_0\). For completeness, we sketch an indirect proof of these facts: Oguiso \cite{Ogu} showed that a non-projective K3 surface \(\tilde{X}\) that admits an automorphism with positive entropy must have its automorphism group map onto \(\mathbb{Z}\) with a finite kernel; taking \(\tilde{X}\) to be the Kummer surface associated to \(X\) and observing that the algebraic automorphism group of \(X\) maps into the automorphism group of \(\tilde{X}\) with a kernel of order two concludes the proof; the precise possibilities for the algebraic automorphism group of \(X\) can be deduced from direct testing of the possible finite-order actions on \(H^*(X,\mathbb{Z})\).

As indicated in Examples 4.1, 4.2, and 4.3, the algebraic automorphism group of a torus of the form \(E \times E\), where \(E\) is an elliptic curve, is isomorphic to \(\xGL_2(K)\), where \(K\) is the ring of surjective algebraic endomorphisms of \(E\); the entropy of any algebraic automorphism of \(E \times E\) is the logarithm of the square of the spectral radius of its image in \(\xGL_2(K)\). (See also \cite{Mum}, \S IV.19.) Example 4.1 shows that \(E \times E\) always exhibits every entropy \(\log(\lambda)\) where \(\lambda\) is a degree-two Salem number which is the square of a quadratic integer; such a \(\lambda\) appears in every real quadratic field. Examples 4.2 and 4.3 show that \(E \times E\) can exhibit even more positive entropies when \(E\) has complex multiplication.

\subsection{The Set of Tori Exhibiting a Fixed Entropy}

Suppose that \(\sigma_1\) and \(\sigma_2\) are algebraic automorphisms with positive entropy of, respectively, two-dimensional complex tori \(X_1\) and \(X_2\), and suppose further that
\[G = \left( \begin{array}{cc}
\gamma_1 & 0 \\
0 & \gamma_2
\end{array} \right)\]
gives the actions of both \(\sigma_1^*\) on \(H^{1,0}(X_1)\) and \(\sigma_2^*\) on \(H^{1,0}(X_2)\) (so that, for one thing, \(\sigma_1\) and \(\sigma_2\) both have entropy \(\log(\max\{|\gamma_1|^2,|\gamma_2^2|\})\)). So there are lattices \(\Lambda_1\) and \(\Lambda_2\) in \(\mathbb{C}^2\) such that \(G(\Lambda_1)=\Lambda_1\), \(G(\Lambda_2)=\Lambda_2\), \(X_1=\mathbb{C}^2/\Lambda_1\), \(X_2=\mathbb{C}^2/\Lambda_2\), and \(\sigma_1\) and \(\sigma_2\) are the quotients of \(G\) by, respectively, \(\Lambda_1\) and \(\Lambda_2\). If there is an (algebraic) isomorphism \(\phi:X_1 \rightarrow X_2\) such that \(\sigma_2 = \phi \circ \sigma_1 \circ \phi^{-1}\), then there is a matrix \(\Phi \in \xGL_2(\mathbb{C})\) such that \(G \circ \Phi = \Phi \circ G\), \(\Phi(\Lambda_1)=\Lambda_2\), and \(\phi\) is the quotient of \(\Phi\); so
\[G|_{\Lambda_2} = \Phi|_{\Lambda_1} \circ G|_{\Lambda_1} \circ \Phi^{-1}|_{\Lambda_2}.\]
Suppose now that \(\gamma_1\) and \(\gamma_2\) are not real. If there is a matrix \(B \in \xGL_4(\mathbb{Z})\) such that \(G|_{\Lambda_2} = B G|_{\Lambda_1} B^{-1}\), then there are matrices \(C\) and \(D\) in \(\xGL_4(\mathbb{R})\) such that both \(C^{-1} \circ B \circ G|_{\Lambda_1} \circ B^{-1} \circ C\) and \(D^{-1} \circ G|_{\Lambda_1} \circ D\) are equal to
\[\left( \begin{array}{cccc}
\Re \gamma_1 & -\Im \gamma_1 & 0 & 0 \\
\Im \gamma_1 & \Re \gamma_1 & 0 & 0 \\
0 & 0 & \Re \gamma_2 & -\Im \gamma_2 \\
0 & 0 & \Im \gamma_2 & \Re \gamma_2 
\end{array} \right) ;\]
it follows that \(C^{-1} \circ B \circ D\) defines a matrix \(\Phi \in \xGL_2(\mathbb{C})\) that commutes with \(G\) and satisfies \(\Phi(\Lambda_1)=\Lambda_2\). So, since each \(\sigma_j^*\) on \(H^1(X_j,\mathbb{Z})\) is given by \((G|_{\Lambda_j})^T\), \(\sigma_1\) and \(\sigma_2\) are the same automorphism (of the same torus) if and only if \(\sigma_1^*\) on \(H^1(X_1,\mathbb{Z})\) and \(\sigma_2^*\) on \(H^1(X_2,\mathbb{Z})\) are conjugate in \(\xGL_4(\mathbb{Z})\); Example 4.1 shows that this statement does not hold when \(\gamma_1\) and \(\gamma_2\) are real. The following result by Latimer and MacDuffee characterizes the sets of \(\xGL_4(\mathbb{Z})\)-conjugacy classes of certain matrices.
\begin{theorem}[\cite{L_M}]
Let \(P(t) \in \mathbb{Z}[t]\) be a monic irreducible polynomial of degree \(r\). Then the set of \(\xGL_r(\mathbb{Z})\)-conjugacy classes of matrices with characteristic polynomial \(P(t)\) is in bijective correspondence with the set of ideal classes in \(\mathbb{Z}[t]/P(t)\).
\end{theorem}
Since the set of ideal classes is finite for any order in a number field, there are only finitely many \(\xGL_r(\mathbb{Z})\)-conjugacy classes of matrices with characteristic polynomial \(P(t)\) in Theorem 5.7. (See, e.g., \cite{C_R}, \S III.20.)

\subsection{Proof of Theorem 1.2}

\emph{Proof of Theorem 1.2:} Let \(S(t)\) be the minimal polynomial for \(\lambda\). Since there are only finitely many monic reciprocal polynomials in \(\mathbb{Z}[t]\) of degree at most four with no roots off the unit circle, there are only finitely many monic and reciprocal degree-six polynomials in \(\mathbb{Z}[t]\) with \(S(t)\) as a factor and four roots on the unit circle. Let
\[Q(t) = t^6 + at^5 + bt^4 + ct^3 + bt^2 + at + 1\]
be a polynomial in \(\mathbb{Z}[t]\) with \(S(t)\) as a factor and four roots on the unit circle such that \(Q(1)=-m^2\) and \(Q(-1)=n^2\) for some integers \(m\) and \(n\); then any polynomial of the form \(t^4+\dots+1 \in \mathbb{Z}[t]\) with the property that the pairwise products of its distinct roots are the roots of \(Q(t)\) must be one of
\[t^4+jt^3-at^2+kt+1 \text{, } t^4-jt^3-at^2-kt+1,\]
\[t^4+kt^3-at^2-jt+1 \text{, or } t^4-kt^3-at^2+jt+1,\]
where \(j=(1/2)(n+m)\) and \(k=(1/2)(n-m)\).

Let \(P(t)\) be a polynomial of the form \(t^4 + \dots + 1 \in \mathbb{Z}[t]\) such that the roots of \(Q(t)\) are the pairwise products of the distinct roots of \(P(t)\). Then \(P(t)\) is reducible if and only if it has a real root, in which case the multiset of roots of \(P(t)\) must be one of
\[\{\sqrt{\lambda},\sqrt{\lambda},\sqrt{\lambda}^{-1},\sqrt{\lambda}^{-1}\} \text{, } \{\sqrt{\lambda},\sqrt{\lambda},-\sqrt{\lambda}^{-1},-\sqrt{\lambda}^{-1}\},\]
\[\{-\sqrt{\lambda},-\sqrt{\lambda},\sqrt{\lambda}^{-1},\sqrt{\lambda}^{-1}\} \text{, or } \{-\sqrt{\lambda},-\sqrt{\lambda},-\sqrt{\lambda}^{-1},-\sqrt{\lambda}^{-1}\}\]
--so that either \(\sqrt{\lambda}+\sqrt{\lambda}^{-1}\) or \(\sqrt{\lambda}-\sqrt{\lambda}^{-1}\) is an integer and therefore either \(\lambda+2+\lambda^{-1}\) or \(\lambda-2+\lambda^{-1}\) is the square of an integer. So, if case 1 does not hold, then it follows from Theorem 5.1 that there are only finitely many \(\xGL_4(\mathbb{Z})\)-conjugacy classes of matrices in \(\xGL_4(\mathbb{Z})\) with characteristic polynomial \(P(t)\); moreover, given such a conjugacy class and a choice of two roots \(\gamma_1\) and \(\gamma_2\) of \(P(t)\) with \(|\gamma_1 \gamma_2|=1\), there is exactly one two-dimensional complex torus \(X\) that admits an algebraic automorphism \(\sigma\) such that \(\sigma^*\) on \(H^1(X,\mathbb{Z})\) is in the conjugacy class and the eigenvalues for \(\sigma^*\) on \(H^{1,0}(X)\) are \(\gamma_1\) and \(\gamma_2\).

If \(\sigma\) is an automorphism of a two-dimensional complex torus \(X\) with entropy \(\log(\lambda)\), then (as in the proof of Theorem 1.1) the characteristic polynomial for \(\sigma^*\) on \(H^2(X,\mathbb{Z})\) must be some \(Q(t)\) as above, and the characteristic polynomial for \(\sigma^*\) on \(H^1(X,\mathbb{Z})\) must be some corresponding \(P(t)\) as above. \(\Box\)

\begin{remark}
We observe that the finite number of distinct tori in case 2 of Theorem 1.2 can be arbitrarily large: for any \(n \in \mathbb{N}\) and every \(k \in \{0,\dots,n\}\),
\[\left( \begin{array}{cc}
0 & -1 \\
1 & 1+2^n\imath
\end{array} \right)\]
gives an automorphism of \(\mathbb{C}/\mathbb{Z}[2^k\imath] \times \mathbb{C}/\mathbb{Z}[2^k\imath]\) whose entropy is the logarithm of the Salem root of
\[t^4-(1+4^n)t^3-2^{2n+1}t^2-(1+4^n)t+1;\]
since \(\mathbb{C}/\mathbb{Z}[2^k\imath]\) has multiplication by \(2^k\imath\) but does not have multiplication by \(2^l\imath\) for any \(l \in \{0,\dots,k-1\}\), each \(k\) gives a distinct torus.
\end{remark}

\begin{remark}
All of the tori in Remark 5.2 are isogenous to one another. In fact, there are no more than 320 isogeny classes represented in case 2 of Theorem 1.2: if
\[G = \left( \begin{array}{cc}
\gamma_1 & 0 \\
0 & \gamma_2
\end{array} \right)\]
satisfies \(|\gamma_1|^2=\lambda\) and \(G(\Lambda_j)=\Lambda_j\) for two lattices \(\Lambda_1\) and \(\Lambda_2\) in \(\mathbb{C}^2\), then \(G|_{\Lambda_1}\) is conjugate to \(G|_{\Lambda_2}\) over \(\mathbb{Q}\); (as above) it follows that there is a matrix \(\Phi \in \xGL_2(\mathbb{C})\) that commutes with \(G\) and satisfies
\[\Phi(\Lambda_1 \otimes \mathbb{Q}) = \Lambda_2 \otimes \mathbb{Q},\]
so that some multiple of \(\Phi\) gives an isogeny from \(\mathbb{C}^2/\Lambda_1\) to \(\mathbb{C}^2/\Lambda_2\); thus the number of distinct isogeny classes exhibiting the entropy \(\log(\lambda)\) is at most the number of pairs \(\{\gamma_1,\gamma_2\}\) of Galois conjugate degree-four algebraic integers with \(|\gamma_1| > |\gamma_2|\) and \(|\gamma_1 \gamma_2|=1\) such that the products of distinct pairs of elements in \(\{\gamma_1,\gamma_2,\overline{\gamma_1},\overline{\gamma_2}\}\) are the roots of some monic and reciprocal degree-six polynomial \(Q(t)\) (as in the proof of Theorem 1.2) that has the minimial polynomial for \(\lambda\) as a factor, has four roots on the unit circle, and satisfies \(Q(1)=-m^2\) and \(Q(-1)=n^2\) for some integers \(m\) and \(n\).
\end{remark}

\bibliographystyle{plain}
\bibliography{ReschkeP-refs-2014.06.27}

\begin{thebibliography}{10}

\bibitem{Ber}
Kenneth Berg.
\newblock Convolution of invariant measures, maximal entropy.
\newblock {\em Math. Syst. Theory}, 3(2):146--150, 1969.

\bibitem{Bea}
M.~J. Bertin, A.~Decomps-Guilloux, M.~Grandet-Hugot, M.~Pathiaux-Delefosse, and
  J.~P. Schreiber.
\newblock {\em Pisot and Salem Numbers}.
\newblock Birkh{\"a}user Verlag, 1992.

\bibitem{B_L}
Christina Birkenhake and Herbert Lange.
\newblock {\em Complex Tori}.
\newblock Birkh{\"a}user, 1999.

\bibitem{Ca1}
Serge Cantat.
\newblock Dynamique des automorphismes des surfaces projectives complexes.
\newblock {\em C. R. Acad. Sci. Paris S{\'e}r. I Math.}, 328(10):901--906,
  1999.

\bibitem{C_R}
Charles Curtis and Irving Reiner.
\newblock {\em Representation Theory of Finite Groups and Associate Algebras}.
\newblock American Mathemical Society, 1966.

\bibitem{Fri}
Shmuel Friedland.
\newblock Entropy of polynomial and rational maps.
\newblock {\em Ann. Math.}, 133(2):359--368, 1991.

\bibitem{Fuj}
Akira Fujiki.
\newblock Finite automorphism groups of complex tori of dimension two.
\newblock {\em Publ. Res. Inst. Math. Sci.}, 24(1):1--97, 1988.

\bibitem{G_V}
{\'E}tienne Ghys and Alberto Verjovsky.
\newblock Locally free holomorphic actions of the complex affine group.
\newblock In {\em Geometric Study of Foliations (Tokyo, 1993)}, pages 201--217.
  World Scientific, 1994.

\bibitem{Gro}
Mikhail Gromov.
\newblock On the entropy of holomorphic maps.
\newblock {\em Enseign. Math.}, 49:217--235, 2003.
\newblock Manuscript, 1977.

\bibitem{Gue}
Vincent Guedj.
\newblock Propri{\'e}t{\'e}s ergodiques des applications rationelles.
\newblock In {\em Panoramas et synth{\`e}ses: Quelques aspects des syst{\`e}mes
  dynamiques polynomiaux}, volume~30, pages 13--95. Soci{\'e}t{\'e}
  Math{\'e}matique de France, 2010.

\bibitem{Kaw}
Shu Kawaguchi.
\newblock Projective surface automorphisms of positive entropy from an
  arithmetic viewpoint.
\newblock {\em Am. J. Math.}, 130(1):159--186, 2008.

\bibitem{Kro}
Leopold Kronecker.
\newblock Zwei s{\"a}tse {\"u}ber {G}leichugnen mit ganzzahligen
  {C}oefficienten.
\newblock {\em J. reine angew. Math.}, 53:173--175, 1857.

\bibitem{L_M}
Claiborne Latimer and Cyrus MacDuffee.
\newblock A correspondence between classes of ideals and classes of matrices.
\newblock {\em Ann. Math.}, 34:313--316, 1933.

\bibitem{Mc1}
Curtis McMullen.
\newblock Dynamics on {K}3 surfaces: {S}alem numbers and {S}iegel disks.
\newblock {\em J. reine angew. Math.}, 545:201--233, 2002.

\bibitem{Mir}
Rick Miranda.
\newblock {\em Algebraic Curves and Riemann Surfaces}.
\newblock American Mathemical Society, 1995.

\bibitem{Mum}
David Mumford.
\newblock {\em Abelian Varieties}.
\newblock Tata Institute of Fundamental Research, 2012.

\bibitem{Ogu}
Keiji Oguiso.
\newblock Bimeromorphic automorphism groups of non-projective hyperk{\"a}hler
  manifolds--a note inspired by {C.T.} {M}c{M}ullen.
\newblock {\em J. Differ. Geom.}, 78(1):163--191, 2008.

\bibitem{Re1}
Paul Reschke.
\newblock Salem numbers and automorphisms of complex surfaces.
\newblock {\em Math. Res. Lett.}, 19(2):475--482, 2012.

\bibitem{Re2}
Paul Reschke.
\newblock Distinguished line bundles for complex surface automorphisms.
\newblock {\em Transform. Groups}, 19(1):225--246, 2014.

\bibitem{S_M}
Tetsuji Shioda and Naoki Mitani.
\newblock Singular abelian surfaces and binary quadratic forms.
\newblock In {\em Classification of Algebraic Varieties and Compact Complex
  Manifolds}, volume 412, pages 259--287. Springer, 1974.

\bibitem{Si2}
Yakov Sinai.
\newblock Flows with finite entropy.
\newblock {\em Dokl. Akad. Nauk SSSR}, 125:1200--1202, 1959.

\bibitem{Si1}
Yakov Sinai.
\newblock The notion of entropy of a dynamical system.
\newblock {\em Dokl. Akad. Nauk SSSR}, 125:768--771, 1959.

\bibitem{V_U}
Gerard van~der Geer and Keiji Ueno.
\newblock Families of abelian surfaces with real multiplication over {H}ilbert
  modular surfaces.
\newblock {\em Nagoya Math. J.}, 88:17--53, 1982.

\bibitem{Yom}
Yosef Yomdin.
\newblock Volume growth and entropy.
\newblock {\em Isr. J. Math.}, 57(3):285--300, 1987.

\bibitem{Zha}
Shou-Wu Zhang.
\newblock Distributions in algebraic dynamics.
\newblock In {\em Surveys in Differential Geometry: Essays in Geometry in
  Memory of S. S. Chern}, volume~10, pages 381--430. International Press, 2006.

\bibitem{Zuc}
Steven Zucker.
\newblock The {H}odge conjecture for cubic fourfolds.
\newblock {\em Compos. Math.}, 34(2):199--209, 1977.

\end{thebibliography}

\end{document}